\documentclass[a4paper,10pt]{article}
\usepackage{graphicx}
\usepackage{amsfonts}
\usepackage{amssymb}

\newtheorem{thm}{Theorem}[section]
\newtheorem{cor}[thm]{Corollary}

\newtheorem{Definition}[thm]{Definition}

\newtheorem{remark}{Remark}
\begin{document}

\title{\bf
Coefficient  Estimates for Certain Subclass of Bi Functions Associated the Horadam Polynomials
 }
\author{Adnan Ghazy Alamoush \\
}
\date{\sl \small Faculty of Science, Taibah University, Saudi Arabia.\\
}
\maketitle

\maketitle

\begin{abstract}
\begin{flushleft}
 In the present article, our goal is  finding estimates on the Taylor-Maclaurin coefficients $|a_{2}|$ and $|a_{3}|$ for a new class of bi-univalent functions defined by means of the Horadam polynomials.
 Fekete-Szeg\"{o} inequalities of functions belonging to this subclass are also founded.

\end{flushleft}

\end{abstract}
\begin{flushleft}
\end{flushleft}

\noindent
{\bf MSC Subj. Class:} Primary 11B39 . 30C45 . 33C45; Secondary 30C50 . 33C05

\section{Introduction and Preliminaries}
A functions of the form $f(z)$ normalized by the following Taylor Maclaurin series

\begin{equation}\label{eg:}
    f(z)=z+\sum^{\infty}_{n=2} a_{n}z^{n}
\end{equation}
which are analytic in the open unit open disk $\mathbb{U}= \{z : z \in \mathbb{C} ,  |z| <1\}$, and belongs to class $A$. Let $\mathbb{S}$ be class of all functions in $A$ which are univalent and normalized by the conditions
$$f(0)=0=f^{'}(0)-1$$ in $\mathbb{U}$. Some of the important and well-investigated subclasses of the univalent function class
$\mathbb{S}$ includes the class $S^{*}(\alpha)(0 \leq \alpha < 1)$ of starlike functions of order $\alpha$ in $\mathbb{U}$ and the class
$K(\alpha)(0 \leq \alpha < 1)$ of convex functions of order  $\alpha$.For two functions $f$ and $g$, analytic in $\mathbb{U}$, we say that the function $f$ is subordinate to $g$ in $\mathbb{U}$,  written as $f(z)\prec g(x),\ \ \ (z\in \mathbb{U})$, provided that there exists an analytic function(that is, Schwarz function) $w(z)$ defined on $\mathbb{U}$  with $$\ \ \  w(0)=0 \ \ and \ \  \ \ \ \ |w(z)|<1\  for\ \  all\ \ z\in \mathbb{U},$$
such that $f(z)=g(w(z))$ for all $z\in \mathbb{U}$.\\

Indeed,it is known that
$$f(z)\prec g(z)\ (z\in \mathbb{U}) \ \Rightarrow\ f(0)=g(0)\ \ and\ \ f( \mathbb{U})\subset g( \mathbb{U}).$$

It is well known that every function
$f\in \mathbb{S}$ has an inverse $f^{-1}$, defined by
$$f^{-1}(f(z))=z\ \ (z \in \mathbb{U}),$$
and
$$f^{-1}(f(w))=w\ \ (|w|<r_{0}(f); r_{0}(f)\geq\frac{1}{4}),$$
where
\begin{equation}\label{eg:}
f^{-1}(w)=w+a_{2}w^{2}+(2a^{2}_{2}-3a_{3})w^{3}-(5a^{3}_{2}-5a_{2}a_{3}+a_{4})w^{4}+...\  .
\end{equation}
A function $f \in A$ is said to be bi-univalent in $\mathbb{U}$ if both $f(z)$ and $f^{-1}(z)$ are univalent in $\mathbb{U}$. For a brief historical account and for several interesting
examples of functions in the class $\Sigma'$. It may be of interest to recall that Lewin (1967) studied
the class of bi-univalent functions and derived that $|a_{2}|<1.51$ (For more details see \cite{Lewin}). Later the papers of Brannan and Taha \cite{Brannan} and Srivastava et al. \cite{Srivastava} and other (see \cite{Alamoush}, \cite{Alamoush1}, \cite{Alamoush2} ) studied the bi-univalent results for may classes.\\
 Moreover, In 1985 Brannan and Taha \cite{Brannan} introduced certain subclasses of the bi-univalent function class $\Sigma'$
similar to the familiar subclasses $S^{*}(\alpha)$ and $K(\alpha)$ of starlike
and convex functions of order $\alpha\ (0\leq\alpha<1)$ in $\mathbb{U}$. respectively (see also \cite{Netanyahu}). Also, they found non-sharp estimates for the initial for the classes $S_{\Sigma'}^{*}(\alpha)$ and $K_{\Sigma'}(\alpha)$ of bi-starlike functions of order $\alpha$ in $\mathbb{U}$ and bi-convex functions of order $\alpha$ in $\mathbb{U}$, corresponding to the function classes $S^{*}(\alpha)$ and $K(\alpha)$.

 Up to now, for the coefficient estimate problem for each of the following Taylor-Maclaurin coefficients $|a_{n}|(n \in \mathbb{N}\setminus\left\{1,2\right\}$ for each $f \in \Sigma'$ given by (1) is  still an open problem.\\

Recently, Horcum and  Kocer \cite{Tugba} considered Horadam polynomials $h_{n}(x)$, which are
given by the following recurrence relation
\begin{equation}\label{eg:}
  h_{n}(x)=pxh_{n-1}(x)+qh_{n-2}(x),\ \ \ (n\in \mathbb{N}\geq2),
\end{equation}
with $h_{1}=a,\ h_{2}=bx,$ and $ h_{3}=pbx^{2}+aq$ where ($a,b,p,q$ are some real constants)  .\\
The characteristic equation of recurrence relation (3) is
\begin{equation}\label{eg:}
t^{2}-pxt-q=0.\end{equation}
This equation has two real roots;
$$\alpha=\frac{px+\sqrt{p^{2}x^{2}+4q}}{2},$$ and $$\beta=\frac{px-\sqrt{p^{2}x^{2}+4q}}{2}.$$
Particular cases of Horadam polynomials sequence are
\begin{itemize}
  \item If $a = b = p = q = 1$, the Fibonacci polynomials sequence is obtained
  \[F_{n}(x)=xF_{n-1}(x)+F_{n-2}(x),\ F_{1}=1,\ F_{2}=x.\]
  \item If $a = 2, b= p = q = 1$, the Lucas polynomials sequence is obtained
    \[L_{n-1}(x)=xL_{n-2}(x)+L_{n-3}(x),\ L_{0}=2,\ L_{1}=x.\]
  \item If $a = q = 1, b= p = 2$, the Pell polynomials sequence is obtained
  \[P_{n}(x)=2xP_{n-1}(x)+P_{n-2}(x),\ p_{1}=1,\ P_{2}=2x.\]
  \item If $a = b = p = 2, q = 1$, the Pell-Lucas polynomials sequence is obtained
    \[Q_{n-1}(x)=2xQ_{n-2}(x)+Q_{n-3}(x),\ Q_{0}=2,\ Q_{1}=2x.\]
  \item If $a = 1, b = p = 2, q = −1$, the Chebyshev polynomials of second kind sequence is obtained
   \[U_{n-1}(x)=2xU_{n-2}(x)+U_{n-3}(x),\ U_{0}=1,\ U_{1}=2x.\]
  \item  If $a = 1, b = p = 2, q = −1$, the Chebyshev polynomials of First kind sequence is obtained
   \[T_{n-1}(x)=2xT_{n-2}(x)+T_{n-3}(x),\ T_{0}=1,\ T_{1}=x.\]
  \item If $x = 1$, The Horadam numbers sequence is obtained
 \[h_{n-1}(1)=ph_{n-2}(1)+qh_{n-3}(1),\ h_{0}=1,\ h_{1}=b.\]
\end{itemize}
For more information associated with these polynomials sequences in [\cite{Horadam}, \cite{Horadam1}, \cite{Koshy}, \cite{Lupas}].\\
These polynomials, the families of orthogonal polynomials and other special polynomials as well as their
generalizations are potentially important in a variety of disciplines in many of sciences, specially in the mathematics, statistics  and  physics.

\begin{remark} \cite{Horadam}
Let $\Omega(x,z)$ be the generating function of the Horadam polynomials
$h_{n}(x)$. Then
\begin{equation}\label{:}
 \Omega(x,z) =\frac{a+(b-ap)xt}{1-pxt-qt^{2}}=\sum^{\infty}_{n=1} h_{n}(x)z^{n-1}.
\end{equation}
\end{remark}
\begin{remark} 
In its special case when $a = 1, b = p = 2, q = −1$ and $x\rightarrow t$, the generating function in Eq. (5) reduces to that of the Chebyshev polynomials $U_{n}(t)$ of the second kind, which is
given explicitly by (see, \cite{Szego})
\[U_{n}(t)=(n+1)\ _{2}F_{1}\left(-n,n+2;\frac{3}{2}; \frac{1-t}{2}\right)=\frac{sin(n+1)\vartheta}{sin(\vartheta)},\ \ t=sin(\vartheta).\]
\end{remark}

In this paper, Our goal is using the Horadam polynomials $h_{n}(x)$ and
the generating function $\Omega(x,z)$ which are given by the recurrence relation (3) and (5), respectively,  to introduce a new subclass of the bi-univalent function class $\Sigma'$. Also,  we provide the initial coefﬁcients and the Fekete-Szeg\"{o} inequality for functions belonging to the class $\Sigma'(x).$

\section{Coefficient Bounds for the Function Class  $\Sigma'(x)$ }
We begin by introducing the function class  $\Sigma'(x)$ by means of the following definitions.

\begin{Definition}
A function $f \in  \Sigma'$ given by (1) is said to be in the class $\Sigma'(x)$,  if the following conditions are satisﬁed:
\begin{equation}\label{eg:}
f^{'}(z)\prec \Omega(x,z)+1-\alpha
\end{equation}
and
\begin{equation}\label{eg:}
g^{'}(w)\prec \Omega(x,w)+1-\alpha
\end{equation}
where the real constants $a$, $b$ and $q$ are as in (3) and $g(w)=f^{-1}(z)$ is given by (2).
\end{Definition}
We first state and prove the following result.
\begin{thm} Let the function $f\in\Sigma'$ given by (1) be in the class $\Sigma'(x)$. Then
\end{thm}
\begin{flushleft}
\begin{equation}\label{eg:}
|a_{2}|\leq \frac{|bx|\sqrt{|bx|}}{\sqrt{\left|bx^{2}(3b-4p)-4aq\right|}}
\end{equation}
\begin{equation}\label{eg:}
|a_{3}|\leq \frac{|bx|}{3}+\frac{(bx)^{2}}{4},
\end{equation}
and for some $\eta\in \mathbb{R}$,
\begin{equation}\label{eg:}
|a_{3}-\eta a^{2}_{2}| \leq   {\begin{array}{c}\left \{ {\begin{array}{c}\frac{|2bx|}{3}\ \ \ \ \ \ \ \ \ \ ,\ \ |\eta-1|\leq 1-\frac{|4(pbx^{2}+qa)|}{3b^{2}x^{2}}\\ \frac{|bx|^{3}|1-\eta|}{3b^{2}x^{2}-4(pbx^{2}+qa)} \ \ \ \ , \ \ |\eta-1|\geq 1-\frac{|4(pbx^{2}+qa)|}{3b^{2}x^{2}}\end{array}} \right \}.
\end{array}}
\end{equation}

\textbf{Proof.}  Let $f\in\Sigma'$  be given by the Taylor-Maclaurin expansion (1). Then, for some analytic
functions $\Psi$ and $\Phi$ such that $\Psi(0)=\Phi(0)=0,$ $|\psi(z)|<1$ and $|\Phi(w)|<1,\ z,w\in \mathbb{U}$ and using Definition 2.1, we can write

\end{flushleft}
\[f^{'}(z) =\omega(x,\Phi(z))+1-\alpha\]
and
\[g^{'}(w)= \omega(x,\psi(w))+1-\alpha\]
or, equivalently,
\begin{equation}\label{eg:}
f^{'}(z)=1+h_{1}(x)-a+h_{2}(x)\Phi(z)+h_{3}(x)[\Phi(z)]^{3}+...
\end{equation}
and

\begin{equation}\label{eg:}
g^{'}(z)=1+h_{1}(x)-a+h_{2}(x)\psi(w)+h_{3}(x)[\psi(w)]^{3}+...\ .
\end{equation}
From (11) and (12), we obtain
\begin{equation}\label{eg:}
f^{'}(z)=1+h_{2}(x)p_{1}z+[h_{2}(x)p_{2}+h_{3}(x)p^{2}_{1}]z^{2}+...
\end{equation}
and

\begin{equation}\label{eg:}
f^{'}(z)=1+h_{2}(x)p_{1}w+[h_{2}(x)q_{2}+h_{3}(x)q^{2}_{1}]w^{2}+...\ .
\end{equation}
Notice that if
\[|\Phi(z)|=|p_{1}z+p_{2}z^{2}+p_{3}z^{3}+...|<1\ \ \ (z\in \mathbb{U})\]
and
\[|\psi(w)|=|q_{1}w+q_{2}w^{2}+q_{3}w^{3}+...|<1\ \ \ (w\in \mathbb{U}),\]
then
$$|p_{i}|\leq1\ \ \ and \ \ \ |q_{i}|\leq1\ \ \ (i\in \mathbb{N}).$$
Thus, upon comparing the corresponding coefficients in (13) and (14), we have
\begin{equation}\label{eg:}
2a_{2}= h_{2}(x)p_{1}
\end{equation}
\begin{equation}\label{eg:}
3a_{3}= h_{2}(x)p_{2}+h_{3}(x)p^{2}_{1}
\end{equation}
\begin{equation}\label{eg:}
-2a_{2}= h_{1}(x)q_{1}
\end{equation}
\begin{equation}\label{eg:}
3\left[2a_{2}^{2}-a_{3}\right]= h_{2}(x)q_{2}+h_{3}(x)q^{2}_{1}.
\end{equation}
From (15) and (17), we find that
\begin{equation}\label{eg:}
p_{1}= -q_{1}
\end{equation}
and
\begin{equation}\label{eg:}
8a^{2}_{2}=h^{2}_{1}(x)(p^{2}_{1}+q^{2}_{1})
\end{equation}

Also, by using (16) and (18), we obtain
\begin{equation}\label{eg:}
6a^{2}_{2}=h_{2}(x)(p_{2}+q_{2})+h_{3}(x)(p^{2}_{1}+q^{2}_{1}).
\end{equation}
By using (19) in (20), we get
\begin{equation}\label{eg:}
\left[6+\frac{8h_{3}(x)}{[h_{2}(x)]^{2}}\right]a^{2}_{2}=h_{2}(x)(p_{2}+q_{2}).
\end{equation}
From (5), (13) and (21), we have the desired inequality (8).
\\

Next, by subtracting (18) from (16), we have
\begin{equation}\label{eg:}
6\left[a_{3}-a^{2}_{2}\right]=h_{2}(x)(p_{2}-q_{2})+h_{3}(x)(p^{2}_{1}-q^{2}_{1}).
\end{equation}
In view of (19) and (20), Equation (23) becomes
\begin{equation}\label{eg:}
a_{3}=a^{2}_{2}+\frac{h_{2}(x)(p_{2}-q_{2})}{6}.
\end{equation}
Hence using (19) and applying (5), we get desired inequality (9).
\\
Now, by using (21) and (23) for some $\eta\in \mathbb{R}$, we get
\[a_{3}-\eta a^{2}_{2}=\frac{[h_{2}(x)]^{3}(1-\eta)(p_{2}+ q_{2})}{6[h_{2}(x)]^{2}-8h_{3}(x)}+\frac{h_{2}(x)(p_{2}- q_{2})}{6}\]
\[=h_{2}(x)\left[\left(\Theta(\eta,x)+\frac{1}{6}\right)p_{2}+\left(\Theta(\eta,x)-\frac{1}{6}\right)q_{2}\right],\]
where
\[\Theta(\eta,x)=\frac{[h_{2}(x)]^{2}(1-\eta)}{6[h_{2}(x)]^{2}-8h_{3}(x)}.\]
So, we conclude that
$$ |a_{3}-\eta a^{2}_{2}| \leq   {\begin{array}{c}\left \{ {\begin{array}{c}\frac{h_{2}(x)}{(3)}\ \ \ \ , \ |\Theta(\eta,x)|\leq \frac{1}{6)}\\ 4h_{2}(x)|\Theta(\eta,x)| \ \ \ \ , \ \ |\Theta(\eta,x)|\geq \frac{1}{6}\end{array}} \right \}
\end{array}} .  $$\\

This proves Theorem 2.2.\\

In light of Remark 2, Theorem 2.2 would yield the following known result.

\begin{cor}
 For  $t\in(1/2,1)$, let the function $f\in\Sigma'$ given by (1) be in the class $\Sigma'(t)$. Then
\begin{equation}\label{eg:}
|a_{1}|\leq \frac{t\sqrt{2t}}{\sqrt{\left|1-t^{2}\right|}}
\end{equation}
\begin{equation}\label{eg:}
|a_{3}|\leq \frac{2t}{3}+t^{2},
\end{equation}
and for some $\eta\in R$,
\begin{equation}\label{eg:}
|a_{3}-\eta a^{2}_{2}| \leq   {\begin{array}{c}\left \{ {\begin{array}{c}\frac{4t}{3}\ \ \ \ , \ |\eta-1|\leq \frac{1-t^{2}}{3t^{2}}\\ \frac{2|\eta-1|}{1-t^{2}} \ \ \ \ , \ \ |\eta-1|\geq \frac{1-t^{2}}{3t^{2}}\end{array}} \right \}.
\end{array}}
\end{equation}

\end{cor}
Taking $\eta=1$ in Corollary 2.3, we get the following consequence

\begin{cor}
 For  t∈(1/2,1), let the function $f\in\Sigma'$ given by (1) be in the class $\Sigma'(t)$. Then
\begin{equation}\label{eg:}
|a_{3}- a^{2}_{2}| \leq \frac{4t}{3}.
\end{equation}

\end{cor}

\ \\

%
%
\end{document}